\documentclass[12pt]{article}
\usepackage{amsmath, latexsym, amssymb}

\numberwithin{equation}{section}

\newcommand{\QED}{\hspace{.2in}\square\newline}

\newtheorem{proposition}{Proposition}[section]
\newtheorem{definition}{Definition}[section]

\newtheorem{claim}{Claim}[section]

\begin{document}

\begin{center}
{\Large \textbf{Exact Summatory Functions for Prime $k$-tuples}} \vskip 4em

{ J. LaChapelle}\\
\vskip 2em
\end{center}

\begin{abstract}
Exact summatory functions that count the number of prime $k$-tuples up to some cut-off integer are presented. Related $k$-tuple analogs of the first and second Chebyshev functions are then defined.
\end{abstract}

\section{Introduction}
Interest in counting prime $k$-tuples for $k\geq2$ is a familiar story. A notable contribution to the storyline came from Hardy and Littlewood in their influential paper \cite{HL} where they conjectured remarkably accurate asymptotics. More contemporary works (see e.g. \cite{GR},\cite{GO},\cite{KO} and references therein) mainly utilize sieve methods to progress, and Yitang Zhang's recent breakthrough \cite{Z} has generated considerable excitement.

In this letter we present exact summatory functions that count the number of prime $k$-tuples up to some cut-off integer. The construction is based on the rather trivial observation that the arithmetic function $\mu(n)\Lambda(n)/\log(n)$ furnishes a characteristic function of primes. Importantly, it can be extended to localize onto prime $k$-tuples. 

This characteristic function of primes allows for a direct representation of the exact prime counting function up to integer $x$
\begin{equation*}
\pi(x)=-\sum_{n\leq x}\mu(n)\frac{\Lambda(n)}{\log(n)}\;.
\end{equation*}
Note this is not a Moebius inversion. The representation is readily extended to prime $k$-tuples. For example the exact twin-prime counting function is
\begin{equation*}
\pi_2(x)=\sum_{n\leq x}\mu(n)\mu(n+2)\frac{\Lambda(n)\Lambda(n+2)}{\log(n)\log(n+2)}\;.
\end{equation*}

In addition to counting, the characteristic function suggests $k$-tuple analogs of the Chebyshev functions. Hopefully, possessing these explicit sums will enable further developments since they are constructed from well-studied objects. As a simple example we give a bound on the prime double Chebyshev functions averaged over prime doubles and show that Zhang's result implies they diverge as $x\rightarrow\infty$.

\section{Counting $k$-tuples}
\begin{proposition}
Let $\mathfrak{P}_k$ be the set of prime $k$-tuples, and denote a prime $k$-tuple by $\mathfrak{P}_k\ni\mathfrak{p}_k=\left(p,p+h_2,\ldots,p+h_k\right)$ with $\mathcal{H}_k:=\{0,h_2,\ldots,h_k\}$ not necessarily admissible. The number of prime $k$-tuples up to some cut-off integer $ x\geq p+h_k$ is given by\footnote{The subscript $(k)$ is supposed to indicate both the order $k$ of the prime tuple and (implicitly) an associated $\mathcal{H}_k:=\{0,h_2,\ldots,h_k\}$. We will make the dependence on $\mathcal{H}_k$ explicit when necessary.}
\begin{eqnarray}
\pi_{(k)}( x)&:=&\sum^{ x}_{\mathfrak{p}_k\in\mathfrak{P}_k}1\notag\\
&=&(-1)^k\sum_{n=2}^{ x}\mu(n)\cdots\mu(n+h_k)\frac{\Lambda(n)}{\log(n)}
\cdots\frac{\Lambda(n+h_k)}{\log(n+h_k)}\;.\notag\\
\end{eqnarray}
In particular, the number of prime doubles $(p,p+2i)$ such that $ x-2\geq2i\in\mathbb{N}_+$ is
\begin{equation}
\pi_{(2)}( x):=\sum^{ x}_{\mathfrak{p}_2\in\mathfrak{P}_2}1=\sum_{n=2}^{ x}\mu(n)\mu(n+2i)\frac{\Lambda(n)}{\log(n)}
\frac{\Lambda(n+2i)}{\log(n+2i)}
\end{equation}
with twin primes corresponding to $i=1$.
\end{proposition}

\emph{Proof}:
Since $\Lambda(n)$ restricts to prime powers $p^\nu$ while $\mu(p^\nu)$ allows only $\nu=1$, then
\begin{equation}
 \mu(n)\Lambda(n)=\left\{\begin{array}{c}
                    -\log(p)\;\;\;\;\;n=p\in\mathfrak{P}_1 \\
                    0\;\;\;\;\;\;\;\;\;\;\;\;\mathrm{otherwise}
                  \end{array}\right.\;.
\end{equation}
Loosely, $\mu(n)\Lambda(n)/\log(n)$ acts like a Dirac delta function for primes relative to the discrete measure on natural numbers. More precisely,
\begin{equation}\label{prime sum}
-\sum_{n=2}^{ x}\mu(n)\,\frac{\Lambda(n)}{\log(n)}=\sum^{ x}_{\mathfrak{p}_1\in\mathfrak{P}_1}1\;.
\end{equation}
Simple induction on $ x$ proves (\ref{prime sum}) since it is obviously true for $ x=2$ and it jumps by one iff $ x+1\in\mathfrak{P}_1$.

In general let $\mathfrak{n}_k:=(n,\ldots,n+h_k)$, then
\begin{eqnarray}
 \mu(n)\Lambda(n)\cdots\mu(n+h_k)\Lambda(n+h_k)=\left\{\begin{array}{l}
(-1)^{k}\log(p)\cdots\log(p+h_k)\,,\;\mathfrak{n}_k=\mathfrak{p}_k\in\mathfrak{P}_{k}\\
0\hspace{2.0in}\mathrm{otherwise}\;.
                         \end{array}\right.\notag\\
\end{eqnarray}
Viewing $\mathfrak{n}_k$ as a point in a $k$-lattice and $\mathfrak{P}_{k+1}$ as a subset of $\mathfrak{P}_{k}\times\mathbb{N}_+=\bigotimes_{k}\mathfrak{P}_{1}\times\mathbb{N}_+$, the $k$-tuple result follows after observing that
\begin{eqnarray}
&&\sum_{n=2}^{ x}\left[\frac{\mu(n)\Lambda(n)\cdots\mu(n+h_k)\Lambda(n+h_k)}{\log(n)\cdots\log(n+h_k)}\right]
\frac{\mu(n+h_{k+1})\Lambda(n+h_{k+1})}{\log(n+h_{k+1})}\notag\\
\hspace{.5in}&&=(-1)^k\sum_{\begin{array}{c}
            \scriptstyle{n'\leq x+h_{k+1}} \\
            \scriptstyle{\mathfrak{n}_k\in\mathfrak{P}_k}
          \end{array}}
\frac{\mu(n')\Lambda(n')}{\log(n')}\,\delta(n'\,,\,(n+h_{k+1}))\notag\\
          \hspace{.5in}&&=(-1)^{k+1}\sum^{ x}_{\mathfrak{p}_{k+1}\in\mathfrak{P}_{k+1}}1\;.
\end{eqnarray}  $\QED$

It is useful to introduce a more compact notation
\begin{equation}
\mu_{(k)}(n):=(-1)^k\mu(n)\cdots\mu(n+h_k)
\end{equation}
and
\begin{equation}
\lambda_{(k)}(n):=\Lambda(n)\cdots\Lambda(n+h_k)/\log(n)\cdots\log(n+h_k)\;.
\end{equation}
So we may write
\begin{equation}
\pi_{(k)}( x)=\sum_{n=2}^{ x}\mu_{(k)}(n)\lambda_{(k)}(n)\;.
\end{equation}

Now define the first and second Chebyshev functions for prime doubles;
\begin{definition}
\begin{eqnarray}
\psi_{(2)}( x)&:=&\frac{1}{2}\sum_{n=2}^{ x}\lambda_{(2)}(n)\log\left(n(n+2i)\right)\;.\\
\theta_{(2)}( x)&:=&\frac{1}{2}\sum_{n=2}^{ x}\mu_{(2)}(n)\lambda_{(2)}(n)\log\left(n(n+2i\right)\;.
\end{eqnarray}
\end{definition}
There are obvious analogs of Chebyshev for higher $k$
\begin{definition}
\begin{eqnarray}
\psi_{(k)}( x)&:=&\sum_{n=2}^{ x}\lambda_{(k)}(n)\log(n_{(k)})
=:\sum_{n=2}^{ x}\frac{\Lambda_{(k)}(n)}{\log^{k-1}(n_{(k)})}\notag\\
\theta_{(k)}( x)&:=&\sum_{n=2}^{ x}\mu_{(k)}(n)\lambda_{(k)}(n)\log(n_{(k)})
=\sum_{n=2}^{ x}\mu_{(k)}(n)\frac{\Lambda_{(k)}(n)}{\log^{k-1}(n_{(k)})}
\end{eqnarray}
where
\begin{equation}
n_{(k)}:=\left(n(n+h_2)\cdots(n+h_k)\right)^{1/k}
\end{equation}
and
\begin{equation}
\Lambda_{(k)}(n):=\lambda_{(k)}(n)\log^{k}(n_{(k)})\;.
\end{equation}
\end{definition}

\begin{proposition}
\begin{equation}
\theta_{(2)}( x)=\frac{1}{2}\sum^{ x}_{\mathfrak{p}_2\in\mathfrak{P}_2}\log\left(p(p+2i)\right)
\end{equation}
\end{proposition}

\emph{Proof}: Use the same reasoning as the previous proof. $\QED$

We can obtain a tight bound on the average (with respect to $i$) prime-double Chebyshev functions. For example,
\begin{eqnarray}
\widehat{\theta_{(2)}}( x)&:=&\frac{\sum_{i=2}^{ x-2}\theta_{(2)}( x)}{\sum_{i=2}^{ x-2}}\notag\\
&=&\frac{1}{2}\frac{1}{( x/2-2)}\sum^{ x}_{\mathfrak{p}_2\in\mathfrak{P}_2}\left[\log(2^{ x/2-2})+\log(p^{ x/2-2})
+\log\left(\frac{\Gamma(\frac{ x+p}{2})}{\Gamma(\frac{4+n}{2})}\right)\right]\notag\\
&\geq&\frac{1}{2}\sum^{ x}_{\mathfrak{p}_2\in\mathfrak{P}_2}\left[\log(2)+\log(3)
+\log\left(\frac{\Gamma(\frac{ x+p}{2})}{\Gamma(\frac{4+p}{2})}\right)^{\frac{1}{ x/2-2}}\right]\notag\\
&>&\frac{1}{2}\widehat{\pi_{(2)}}( x)+\frac{\widehat{\pi_{(2)}}( x)}{ x}\left[\log(\Gamma( x/2+1)-1\right]\notag\\
&=&\widehat{\pi_{(2)}}( x)\left[(O(\log( x))+O(1)\right]\;.
\end{eqnarray}
On the other hand,
\begin{eqnarray}
\frac{1}{2}\sum^{ x}_{\mathfrak{p}_2\in\mathfrak{P}_2}\log\left(p(p+2i)\right)
<\sum^{ x}_{\mathfrak{p}_2\in\mathfrak{P}_2}\log\left(p+2i\right)
\leq\sum^{ x}_{\mathfrak{p}_2\in\mathfrak{P}_2}\log( x)
&=&\log( x)\,\pi_{(2)}( x)\;.\notag\\
\end{eqnarray}
So $\widehat{\theta_{(2)}}( x)\asymp\log( x)\,\widehat{\pi_{(2)}}( x)$. Because of Zhang's theorem \cite{Z}, $\widehat{\pi_{(2)}}( x)$ must diverge with $ x$. It follows that $\lim_{ x\rightarrow\infty}\widehat{\theta_{(2)}}( x)/\log( x)=\infty$.  Clearly the same bounds obtain for $\widehat{\psi_{(2)}}( x)$ in terms of $\widehat{J_{(2)}}( x)$ where $J_{(2)}$ is the weighted sum of prime-power doubles.

If the gamma hypothesis \cite{LA1} holds, this can be strengthened:
\begin{claim}\label{claim}
Assume the gamma distribution hypothesis  and that $\mathcal{H}_k$ is admissible. Then, for each $k\in\mathbb{N}_+$, $\pi_{(k)}( x)\sim c_{(k)}\left( x/\log^k( x)\right)$ for some positive constant  $c_{(k)}$.
\end{claim}

\emph{Sketch of proof}\footnote{The idea of the proof is obvious and simple since the gamma distribution hypothesis eliminates the difficult part of the proof by fiat. Of course strict analytic rigor is necessary to promote this claim to a theorem: even then it would hold only conditionally.}:
As a direct consequence of the PNT;
\begin{equation}
\begin{array}{c}
  \pi_{(1)}( x)=-\sum^{ x}_{n=2}\mu_{(1)}(n)\lambda_{(1)}(n)\sim x/\log( x)
   \vspace{.05in}\\
  \big\Downarrow \\
  \pi_{(1)}(N+1)-\pi_{(1)}(N)\sim 1/\log(N)\;.
\end{array}
\end{equation}
Hence the density of log-primes goes like $\mu(n)\Lambda(n)\sim1$.

Now, relative to the density of prime doubles, $\pi_{(2)}( x)$ represents a weighted intersection between  the $2$-lattice of integers and a ray $\mathfrak{r}_{(2i)}$ along a direction dictated by $\mathcal{H}_2=\{0,2i\}$. Typically, $\mathfrak{r}_{(2i)}$ will meet only coprime points due to the Von Mangoldt product so it is more efficient to restrict attention to the coprime $2$-lattice.  Explicitly, denote the coprime $2$-lattice by the set of points $\{(n_1,n_2)\in\mathbb{N}_+^2 \;|\; \mathrm{gcd}(n_1,n_2)=1\}$. Then
\begin{eqnarray}
\pi_{(2i)}( x)&=&\sum^{ x}_{n=2}\mu_{(2i)}(n)\lambda_{(2i)}(n)\notag\\
&=&\sum_{n_1\leq x}\sum_{n_2\leq x+2i}\mu_{(1)}(n_1)\lambda_{(1)}(n_1)\mu_{(1)}(n_2)\lambda_{(1)}(n_2)
\,\delta(n_2\,,\,n_1+2i)\notag\\
&=&\sum_{n_1\leq x}\frac{\mu(n_1)\Lambda(n_1)}{\log(n_1)}\,
\left(\sum_{n_2\leq x+2i}\frac{c_{(2i)}(n_2|n_1)}{\log(n_2)}
\,\delta(n_2\,,\,n_1+2i)\right)\;.
\end{eqnarray}
The delta function restricts to $\mathfrak{r}_{(2i)}$ where $n_1$ and $n_1+2i$ are coprime. The weight function $c_{(2)}(n_2|n_1)$ encodes the condition $\mu(n_2)\Lambda(n_2)\neq0$ given that $\mu(n_1)\Lambda(n_1)\neq0$ along $\mathfrak{r}_{(2i)}$.

Consider a change of summation variable $n_2\mapsto n_2'=n_1+2j$ in the inner sum with $2j=(2i)^l$ ensuring that $n'_2$ is also coprime to $n_1$. This will generate a different ray $\mathfrak{r}_{(2j)}$, but the intersections of $\mathfrak{r}_{(2i)}$ and $\mathfrak{r}_{(2j)}$ with the lattice are congruent in the sense that both intersect the lattice for the same set of $n_1$ because $\mathrm{gcd}(2i,2j)=2i$. In essence, $\mathfrak{r}_{(2i)}$ gets shifted `up' $2j-2i$ units. But, according to the gamma hypothesis, the counting of prime powers is a random process. The two counting processes along $\mathfrak{r}_{(2i)}$ and  $\mathfrak{r}_{(2j)}$ are independent, and they share congruent sets of coprime events and the same underlying distribution of prime powers. Consequently, for sufficiently large $ x$, we expect $c_{(2i)}(n_2|n_1)$ to depend only on the equivalence class $[2i]$ defined by the relation $\mathrm{gcd}(2i,2j)=2i$.
Hence, for sufficiently large $N$,
\begin{equation}
\pi_{[2i]}(N+1)-\pi_{[2i]}(N)\sim\frac{1}{\log(N)}\frac{c_{[2i]}(N)}{\log(N)}\;.
\end{equation}

Note that $c_{[2i]}(N)$ applies to counting numbers along $\mathfrak{r}_{[2i]}$, and it detects prime doubles. It is decisive here that the assumed gamma distribution for prime doubles \cite{LA2} allows this condition on counting numbers to be quantified as a probability statement implying that $c_{[2i]}(N)=O(1)$ for all $i$.\footnote{Since we have assumed a probability model, the usual identification of $c_{[2i]}(N)$ with the singular series can be made.} In other words, every $\mathfrak{r}_{(2)}$ has a non-vanishing  intersection with the coprime $2$-lattice for \emph{any} sufficiently large region of the lattice, and the intersection density is asymptotically constant. So we end up with $\pi_{(2)}( x)\sim c_{(2)}\left( x/\log^2( x)\right)$. The argument easily extends inductively, so $\pi_{(k)}( x)\sim c_{(k)}\left( x/\log^k( x)\right)$ in general. $\QED$

Aside from a few details, this sketch is essentially the familiar probabilistic argument so it doesn't carry much weight. However, the exercise is useful because it is formulated in the context of a pair-wise coprime lattice structure which affords an unconventional perspective.

\end{document}